\documentclass[psamsfonts]{amsart}

\usepackage{amssymb,amsfonts}
\usepackage[all,arc]{xy}
\usepackage{enumerate}
\usepackage{mathrsfs}

\newtheorem{theorem}{Theorem}[section]
\newtheorem{corollary}[theorem]{Corollary}
\newtheorem{proposition}[theorem]{Proposition}

\newtheorem{questions}[theorem]{Questions}

\theoremstyle{definition}
\newtheorem{definition}[theorem]{Definition}

\theoremstyle{remark}
\newtheorem{remark}[theorem]{Remark}

\let\phi=\varphi

\def\card{\operatorname{card}}

\let\oldbigwedge\bigwedge
\def\BIGwedge{{\textstyle\oldbigwedge}}
\def\medwedge{{\scriptstyle\oldbigwedge}}
\def\bigwedge{\mathchoice{\BIGwedge}{\BIGwedge}{\medwedge}{}}

\DeclareMathOperator{\Av}{Av}

\DeclareMathOperator{\f}{frac}

\makeatletter

\let\epsilon=\varepsilon

\let\c@equation\c@theorem
\makeatother
\numberwithin{equation}{section}

\begin{document}
	\title{A Computational Criterion for the Irrationality of Some Real Numbers}
	
	\author[P. Nasehpour]{Peyman Nasehpour}
	
	\address{Department of Engineering Science \\ Golpayegan University of Technology  \\
		Golpayegan\\
		Iran}
	\email{nasehpour@gut.ac.ir, nasehpour@gmail.com}
	
	
	
	\begin{abstract}
		In this paper, we compute the asymptotic average of the decimals of some real numbers. With the help of this computation, we prove that if a real number cannot be represented as a finite decimal and the asymptotic average of its decimals is zero, then it is irrational. We also show that the asymptotic average of the decimals of simply normal numbers is 9/2.
	\end{abstract}
	
	\maketitle
	
	\section{Introduction}
	
	One of the interesting topics in the theory of real numbers is to decide the irrationality of a real number based on the properties of the sequence of its decimal expansions \cite{Bundschuh1984,CsapodiHegyvari2018,Hegyvari1993,Mahler1981,Martinez2001,Mercer1994,Sander1995}. We mention two of such beautiful results:
	
	\begin{itemize}
		
		\item Bundschuh in \cite{Bundschuh1984}, as a generalization of Mahler's theorem in \cite{Mahler1981}, proves that if $g,h\geq 2$ are two fixed integers, then the positive real number $0.(g^0) (g^1) (g^2) \ldots (g^n) \ldots$ is irrational, where $(g^n)$ is to mean the number $g^n$ written in base $h$. 
		
		\item Hegyv\'{a}ri in \cite{Hegyvari1993} proves that if $(a_n)_{n\in \mathbb N}$ is a strictly increasing sequence of positive integers for which $\sum_{n=1}^{\infty} 1/a_n = \infty$, then the decimal fraction $\alpha = 0.(a_1) (a_2) \ldots (a_n) \ldots$ is irrational, where $(a_n)$ is to mean the positive integer number in base 10.
		
	\end{itemize}
	
	Hardy and Wright (cf. \cite[Theorem 137]{HardyWright1975}) show that the real number \[r = 0.r_1r_2r_3 \ldots r_n \ldots = 0.011010100010 \ldots,\] is irrational, where $r_n = 1$ if $n$ is prime and $r_n = 0$ otherwise. Their beautiful proof is based on this fact that no non-constant polynomial in $\mathbb Z[X]$ is a prime-representing function \cite[Theorem 21]{HardyWright1975}. We recall that a function $f(x)$ is said to be a prime-representing function if $f(x)$ is a prime number for all positive integral values of $x$ \cite{Mills1947}. In Section \ref{sec:criterion}, by The Chebychev's Estimate Theorem \cite[Theorem 4.2.1]{FineRosenberger2007} and calculating the asymptotic average of the decimals of some real numbers, we give an alternative proof for Hardy and Wright's theorem (see Theorem \ref{HardyWright}). Note that this is a corollary of the main theorem of this paper which says that if a real number $r$ cannot be represented as a finite decimal and the asymptotic average of its decimals is zero, then $r$ is irrational (see Theorem \ref{AverageDecimalsThm2}).
	
	In Corollary \ref{AverageDecimalsCor}, we show that if $(a_n)_{n\in \mathbb N}$ is a strictly increasing sequence of positive integers such that $\displaystyle \lim_{n\to+\infty}\frac{n}{a_n} =0$, then  $\displaystyle r=\sum_{n=1}^{+\infty}  \frac{b_n}{10^{a_n}}$ is irrational, where each $1 \leq b_n \leq 9$ is a positive integer.
	
	Section \ref{sec:normal} is devoted to the asymptotic average of the decimals of simply normal numbers. Let us recall that a real number $r$ is a simply normal number to base $b$ if for the decimals $(r_n)_{n\in \mathbb N}$ of the fractional part $(0.r_1r_2r_3 \dots r_n \dots)_b $ of the real number $r$, we have the following property:
	
	\[\lim_{n\to+\infty} \frac{\card \{j : 1\leq j \leq n, r_j=d\}}{n} = \frac{1}{b}, \] where $d \in \{0,1,2,\dots,b-1\}$ \cite[Definition 4.1]{Bugeaud2012}. In Theorem \ref{SimplyNormalCor}, we prove that if $r$ is a simply normal number to base $b$, then the asymptotic average of the decimals of $r$ is equal to $\displaystyle \frac{b-1}{2}$. We need to recall that a fractional part of a real number $r$ is the non-negative real number $\f(r) :=|r| - \lfloor |r| \rfloor$, where $|r|$ is the absolute value of $r$ and $\lfloor r \rfloor$ is the integer part of $r$. 
	
	In Proposition \ref{ChampernowneLiouville}, we show that the asymptotic average of the decimals of Champernowne number is 9/2. Note that Champernowne number is the number $C_{10} =0.1234567891011121314151617181920\dots,$ whose sequence of decimals is the increasing sequence of all positive integers. Since we are not aware of this point if the asymptotic average of the decimals of all irrational numbers exists and if it exists we do not know of a systematic method to calculate it, we propose a couple of questions at the end of the paper (check Questions \ref{AsymptoticAveragesQ}).
	
	\section{A Criterion for the Irrationality of Some Real Numbers}\label{sec:criterion}
	
	First we recall some facts related to real numbers in order to fix some definitions and terminologies. Let us recall that every real number to base $b$ can be expressed by a decimal expansion, and this expansion can be performed in only one way \cite[p. 38]{StewartTall2015}. To be more precise, a real number is regular with respect to some base number $b$ when it can be expanded in the corresponding number system with a finite number of negative powers of $b$ \cite[p. 316]{Ore1948}. A regular number is also called a real number with finite decimal \cite[p. 25]{Havil2003}. It is easy to see that a real number $r$ with respect to some base number $b$ is regular if and only if there are coprime integer numbers $p$ and $q$ such that $\displaystyle r = \frac{p}{q}$ and $q$ contains no other prime factors than those that divide $b$ \cite[p. 316]{Ore1948}. In this paper, if the fractional part of a regular number to base $b$ is \[\alpha = (0.r_1 r_2 \dots r_n)_b,\] we only consider the representation with an infinite series of $(b-1)$: \[\alpha = (0.r_1 r_2 \dots (r_n - 1) \overline{(b-1)} \dots)_b .\] Therefore, if we agree always to pick the non-terminating expansion in the case of regular numbers, then fractional part of each real number to base $b$ corresponds uniquely to an infinite decimal $(0.r_1r_2r_3\dots r_n\dots)_b$. Finally, we assert that if $(r_n)_{n\in \mathbb N}$ is a sequence in real numbers, the sequence of the averages is defined as follows: \[a_n = \frac{r_1 + r_2 + \dots + r_n}{n}.\]
	
	\begin{definition}
		
		\label{AverageDecimalsDef}
		
		Let the fractional part of a real number $r$ to base $b$ be \[(0.r_1 r_2 \dots r_n \dots)_b.\] Then, we define the asymptotic average of the decimals of the real number $r$ by \[\Av_b(r) = \lim_{n\to+\infty} \frac{r_1 + r_2 + \dots + r_n}{n},\] if it exists. Usually, we denote $\Av_{10}(r)$ by $\Av(r)$ if there is no fear of any ambiguity.

	\end{definition}
	
	In the following, we calculate the asymptotic average of the decimals of some real numbers:
	
	\begin{theorem}
		
		\label{AverageDecimalsThm1}
		
		Let the decimals $(r_n)_{n\in \mathbb N}$ of the fractional part \[(0.r_1r_2r_3 \dots r_n \dots)_b \] of a real number $r$ satisfy the following:
		
		\[\lim_{n\to+\infty} \frac{\card \{j : 1\leq j \leq n, r_j=d\}}{n} = \omega_d, \] where $d \in \{0,1,2,\dots,(b-1)\}$, $0 \leq \omega_d \leq 1$, and $\displaystyle \sum^{b-1}_{d=0} \omega_d = 1$. Then, \[\Av_b(r) = \sum^{b-1}_{d=1} (d \cdot \omega_d) .\]
		
		\begin{proof}
			Let $A(d,n) = \{j : 1\leq j \leq n, r_j=d\}$. By assumption, \[\lim_{n\to+\infty} \frac{\card A(d,n)}{n} = \omega_d. \] This means that for any $\epsilon > 0$, there is a natural number $N_d$ such that if $n > N_d$, then \[\left| \frac{\card A(d,n)}{n} - \omega_d\right| < \epsilon.\] Now, if we define $N = \max\{N_d\}^{b-1}_{d=0}$, for each $n > N$, we have the following:\begin{equation} \label{inequality} \displaystyle \omega_d - \epsilon < \frac{\card A(d,n)}{n} < \omega_d + \epsilon.  \end{equation}
			
			Since \[ \frac{r_1 + r_2 + \dots + r_n}{n} =\displaystyle \frac{\sum^{b-1}_{d=0} \sum_{i\in A(d,n)} r_i}{n} =  \frac{\sum^{b-1}_{d=1} d\card A(d,n)}{n},\] by using the inequality (\ref{inequality}), we have the following:
			\[  \left| \frac{r_1 + r_2 + \dots + r_n}{n} - \sum^{b-1}_{d=1} d \cdot \omega_d \right| < \frac{b(b-1)\epsilon}{2}. \] Hence, \[\Av_b(r) = \sum^{b-1}_{d=1} d\cdot \omega_d \] and the proof is complete.
		\end{proof}
	\end{theorem}
	
	Let us recall that any rational number is expressible as a finite decimal (if it is regular) or an infinite periodic decimal. Conversely, any decimal expansion which is either finite or infinite periodic is equal to some rational number \cite[p. 32]{Niven1961}. Since in this paper, we only consider the infinite decimal representation of a finite decimal number, we have the following:
	
	\begin{corollary}
		
		\label{AverageDecimalsLem}
		
		Let $r$ be a rational number. Then $\Av_b(r)$ exists and is a positive rational number. Moreover, if the fractional part of the rational number $r$ to base $b$ is \[(0.r_1r_2 \dots r_n \overline{p_1 p_2 \dots p_m})_b,\] then \[\Av_b(r) = \frac{p_1 + p_2 + \dots + p_m}{m}.\] In particular, if $r$ is regular to base $b$, then \[\Av_b(r) = b-1.\]
	\end{corollary}
	
	Now we prove the main theorem of our paper:
	
	\begin{theorem}[A Criterion for the Irrationality of Some Real Numbers]
		
		\label{AverageDecimalsThm2}
		
		Let $r$ be a real number such that $\Av_b(r) =0$. Then $r$ is irrational.
		
		\begin{proof} From Definition \ref{AverageDecimalsDef}, it is clear that $\Av_b(r)$ is non-negative. Also by Corollary \ref{AverageDecimalsLem}, $\Av_b(r)$ is positive if $r$ is rational. Therefore, if $\Av_b(r) = 0$, then $r$ is irrational and the proof is complete. \end{proof}
	\end{theorem}
	
	\begin{corollary}
		
		\label{AverageDecimalsCor}
		
		Let $(a_n)_{n\in \mathbb N}$ be a strictly increasing sequence of positive integers such that \[\lim_{n\to+\infty}\frac{n}{a_n} =0.\] Define $\displaystyle r=\sum_{n=1}^{+\infty}  \frac{b_n}{10^{a_n}},$ where $1 \leq b_n \leq 9$ is a positive integer, for each $n\in \mathbb N$. Then, $r$ is irrational.
		
		\begin{proof} $\displaystyle \Av(r) \leq 9\cdot \lim_{n\to+\infty}\frac{n}{a_n} =0$. \end{proof}  
	\end{corollary}
	
	\begin{remark}
		By Corollary \ref{AverageDecimalsCor}, the Liouville's constant \[\displaystyle \ell=\sum_{n=1}^{+\infty}  \frac{1}{10^{n!}},\] \cite[Theorem 6.6]{Stark1987} is irrational. 
	\end{remark}
	
	Let us recall that $\pi(n) = \card\{p \in \mathbb P: p\leq n\}$ is the prime counting function, where $\mathbb P$ is the set of all prime numbers. The Chebychev's Estimate Theorem states that there exist positive constants $A_1$ and $A_2$ such that \[A_1 \cdot \frac{n}{\ln n} < \pi(n) < A_2 \cdot \frac{n}{\ln n}, \] for all $n \geq 2$ \cite[Theorem 4.2.1]{FineRosenberger2007}. One of the nice corollaries of this theorem says that $\lim_{n\to+\infty} \displaystyle \frac{\pi(n)}{n} = 0$ \cite[Corollary 4.2.3]{FineRosenberger2007}. We use this to give an alternative proof for the following result brought in the book by Hardy and Wright \cite{HardyWright1975}:
	
	\begin{theorem}\label{HardyWright} \cite[Theorem 137]{HardyWright1975} The real number \[r = 0.r_1r_2r_3 \dots r_n \dots = 0.011010100010\dots,\] where $r_n = 1$ if $n$ is prime and $r_n = 0$ otherwise, is irrational.
		
		\begin{proof} It is clear that $\displaystyle \frac{r_1 + r_2 + r_3 + \dots + r_n}{n} = \frac{\pi(n)}{n}.$ Since $\displaystyle \Av(r) = \lim_{n\to+\infty} \frac{\pi(n)}{n} = 0$ \cite[Corollary 4.2.3]{FineRosenberger2007}, $r$ is irrational (Corollary \ref{AverageDecimalsCor}). This finishes the proof.
		\end{proof}
	\end{theorem}
	
	\section{The Asymptotic Average of the Decimals of Simply Normal Numbers}\label{sec:normal}

	Let us recall that a real number $r$ is a simply normal number to base $b$ if for the decimals $(r_n)_{n\in \mathbb N}$ of the fractional part $(0.r_1r_2r_3 \dots r_n \dots)_b $ of the real number $r$, we have the following property:
	
	\[\lim_{n\to+\infty} \frac{\card \{j : 1\leq j \leq n, r_j=d\}}{n} = \frac{1}{b}, \] where $d \in \{0,1,2,\dots,b-1\}$ \cite[Definition 4.1]{Bugeaud2012}.
	
	\begin{theorem}
		
		\label{SimplyNormalCor}
		
		Let $r$ be a simply normal number to base $b$. Then $\displaystyle \Av_b(r) = \frac{b-1}{2}$.
		
		\begin{proof} By Theorem \ref{AverageDecimalsThm1}, $\displaystyle \Av_b(r) = \sum^{b-1}_{d=1} (d \cdot \frac{1}{b}) = \frac{b-1}{2}$. \end{proof}
		
	\end{theorem}
	
	Let us recall that a real number $r$ is algebraic if it is the root of a polynomial $f(X) \in \mathbb Z[X]$, otherwise it is transcendental \cite{ErdosDudley1983,Niven1956}. The set of transcendental numbers is uncountable \cite{Cantor1874}. For a masterful exposition of some central results on irrational and transcendental numbers, refer to \cite{Niven1956}.

	\begin{proposition}
		
		\label{ChampernowneLiouville}
		
		The following statements hold:
		
		\begin{enumerate}
			
			\item There is a transcendental number $r$ such that $\Av(r) \neq 0$.
			
			\item There is a transcendental number $r$ such that $\Av(r) = 0$.
			
		\end{enumerate}
		
		\begin{proof}

			$(1)$: The Champernowne number, \[C_{10} =0.1234567891011121314151617181920\dots,\] whose sequence of decimals is the increasing sequence of all positive integers, is a simply normal number (cf. \cite{Champernowne1933} and \cite[Theorem 4.2]{Bugeaud2012}). So, by Theorem \ref{SimplyNormalCor}, $\displaystyle \Av(r) \neq 0$. On the other hand, $C_{10}$ is transcendental \cite{Mahler1937}.
			
			$(2)$: The Liouville's constant $\displaystyle \ell=\sum_{n=1}^{+\infty}  \frac{1}{10^{n!}}$ is transcendental \cite[Theorem 6.6]{Stark1987} while by Corollary \ref{AverageDecimalsCor}, we have $\Av(\ell) = 0$.\end{proof}
		
	\end{proposition}
	
	\begin{remark}
		With the help of Corollary \ref{AverageDecimalsLem}, it is easy to see that if $s\in (0,9] \cap \mathbb Q$, then there is a rational number $r$ such that $\Av(r) = s$. Based on this, the following questions arise:
	\end{remark}
	
	\begin{questions}
		
		\label{AsymptoticAveragesQ}
		
		\begin{enumerate}
			\item Is there any irrational (transcendental) number $r$ such that $\Av(r)$ does not exist?
			
			\item Is there any irrational (transcendental) number $r$ such that $\Av(r)=s$, for an arbitrary $s$ in $(0,9] \cap (\mathbb R - \mathbb Q)$,?
			
			\item Does the asymptotic avergage of the decimals of the irrational number $\sqrt{2}$ exist and if it exists, what is that? The same question arises for other celebrated irrational numbers such as $e$, $\pi$, $\gamma$, and $\log^3_2$, where \[\gamma = \lim_{n\to+\infty} \big( -\ln n + \sum^n_{k=1} \frac{1}{k}\big)\] is the Euler-Mascheroni constant?
		\end{enumerate}
	\end{questions}
	
	\subsection*{Acknowledgments} This work is supported by the Golpayegan University of Technology. Our special thanks go to the Department of Engineering Science at the Golpayegan University of Technology for providing all the necessary facilities available to us for successfully conducting this research.

\end{document}